# Simultaneous Optimization of Launch Vehicle Stage and Trajectory Considering Operational Safety Constraints

Jaeyoul Ko[1], Jaewoo Kim[2], Jimin Choi[3], and Jaemyung Ahn[4]

*Korea Advanced Institute of Science and Technology, Daejeon 34141, South Korea*

**A conceptual design of a launch vehicle involves the optimization of trajectory and stages considering its launch operations. This process encompasses various disciplines, such as structural design, aerodynamics, propulsion systems, flight control, and stage sizing. Traditional approaches used for the conceptual design of a launch vehicle conduct the stage and trajectory designs sequentially, often leading to high computational complexity and suboptimal results. This paper presents an optimization framework that addresses both trajectory optimization and staging in an integrated way. The proposed framework aims to maximize the payload-to-liftoff mass ratio while satisfying the constraints required for safe launch operations (e.g., the impact points of burnt stages and fairing). A case study demonstrates the advantage of the proposed framework compared to the traditional sequential optimization approach.**

## Nomenclature

| | | |
|---|---|---|
| $v_{ex}$ | = | Exhaust velocity, *m/s* |
| $\varepsilon$ | = | Structural mass fraction, - |
| $m$ | = | Launch vehicle mass, *kg* |
| $m_{s,k}$ | = | Structural mass of stage $k$, *kg* |
| $m_{p,k}$ | = | Propellant mass of stage $k$, *kg* |

---

| | | |
|---|---|---|
| $m_{Lift\text{-}off}$ | = | Lift-off mass of launch vehicle, $kg$ |
| $\mu_k$ | = | Initial-to-final mass raio of stage $k$, - |
| $\dot{m}$ | = | Propellant mass flow rate, $kg/s$ |
| $\boldsymbol{r}$ | = | Position vector, $m$ |
| $\boldsymbol{v}$ | = | Velocity vector, $m/s$ |
| $\Delta v$ | = | Velocity increment created by launch vehicle, $m/s$ |
| $\Delta v_{req}$ | = | Velocity increment required for mission, $m/s$ |
| $\Delta v_{Loss}$ | = | Velocity loss, $m/s$ |
| $\lambda$ | = | Latitude, $rad$ |
| $\varphi$ | = | Longitude, $rad$ |
| $h$ | = | Altitude, $m$ |
| $A_e$ | = | Nozzle exit area, $m^2$ |
| $p_a$ | = | Atmospheric pressure, $Pa$ |
| $\rho_a$ | = | Atmospheric density, $kg/m^3$ |
| $C_D$ | = | Drag coefficient, - |
| $S_{ref}$ | = | Reference area, $m^2$ |
| $J_2$ | = | J2 Coefficient of the Earth |
| $R_E$ | = | Equatorial radius of the Earth, $m$ |
| $G.M.$ | = | Gravitational parameter of the Earth, $m^3/s^2$ |
| $\gamma$ | = | Flight path angle, $rad$ |
| $\delta$ | = | Thrust vector control (TVC) angle, $rad$ |
| $\theta$ | = | Pitch angle, $rad$ |
| $\psi$ | = | Yaw angle, $rad$ |
| $i$ | = | Inclination angle, $rad$ |

**Subscript**

| | | |
|---|---|---|
| PL | = | Payload |
| $k$ | = | $k^{th}$ stage of the launch vehicle |

$i / f$        =   Initial/final

**Supersciprt**

\*           =   Values associated with optimal solution

$p$         =   Values associated with phase $p$

# I. Introduction

The performance of launch vehicles plays a critical role in space exploration. Optimization of launch vehicle performance considers multiple disciplines such as structural design, aerodynamics, propulsion, flight control, and stage sizing. During the conceptual design phase, crucial decisions regarding various factors are made, significantly influencing the final design of the launch vehicle. In the conceptual design process of a launch vehicle, stage sizing and trajectory design are two interconnected key elements. Traditional approaches adopt a sequential stage and trajectory optimization procedure, followed by iterations for convergence [1-2]. The designer first allocates the stage masses based on the target orbit, engine characteristics, and structural ratio based on the assumption of velocity loss terms. Then the trajectory optimization using the stage information is conducted considering operational constraints, yielding the computed velocity loss terms. This procedure continues until the assumed and computed velocity loss terms converge. Past studies demonstrated the effectiveness of the approach above based on sequential design and iterations. However, this approach often fails to capture the inter-dependencies between each element, leading to a suboptimal design. Moreover, the computational complexity of these tasks makes the conceptual design procedure challenging.

This paper discusses a framework that can simultaneously optimize the stages and trajectory of a launch vehicle to overcome the challenge. By doing so, one can effectively explore the design space constrained by the launch vehicle's trajectory, leading to improved performance and computational efficiency. As the contribution of this paper, we proposed a simultaneous optimization framework considering the launch vehicle's operational constraints. We also conducted a case study to design a launcher's stage size and trajectory with realistic vehicle data. Performance comparison of the proposed method and the traditional sequential optimization-based approach demonstrated its effectiveness.

The remainder of this paper is organized as follows. Section 2 provides a brief overview of previous studies. Section 3 introduces the integrated framework for optimizing the stage size and launch trajectory. Section 4 presents several case studies to validate and compare the developed framework to traditional approaches. Lastly, Section 5 concludes the paper and discusses potential future research directions.

## II. Preliminaries

### A. Optimal Staging of a Launch Vehicle

The optimal stage sizing of a launch vehicle determines the mass distributions (propellant and stage masses) for each vehicle stage to put the payload into the target orbit. The propulsion system capability (specific impulse or equivalently effective exhaust velocity) and the structural characteristics (mass fraction) are external parameters for the optimal sizing problem. Malina proposed an efficient optimal sizing approach for multi-stage rockets [3], and many studies utilized the approach for the optimal staging of launch vehicles [4-9]. This section provides a concise overview of previous studies on optimal stage sizing.

We assume the exhaust velocity of the propulsion system and the structural mass ratio for each stage are given parameters used in the optimal stage sizing problem. Note that the structural mass fraction for stage $k$ ($\varepsilon_k$) is the ratio of the structural mass ($m_{s,k}$) to the sum of the structural and propellant ($m_{p,k}$) masses defined as follows.

$$\varepsilon_k \equiv \frac{m_{s,k}}{m_{s,k} + m_{p,k}} \tag{1}$$

Additionally, we define $m_{pl,k}$ as the mass considered as the payload of the $k^{th}$ stage, equivalent to the sum of the masses of the upper stages and payload. Note that the value corresponds to the initial mass of the $(k+1)^{th}$ stage of the launch vehicle ($m_{i,k+1}$), which is the mass immediately after separating the $k^{th}$ stage.

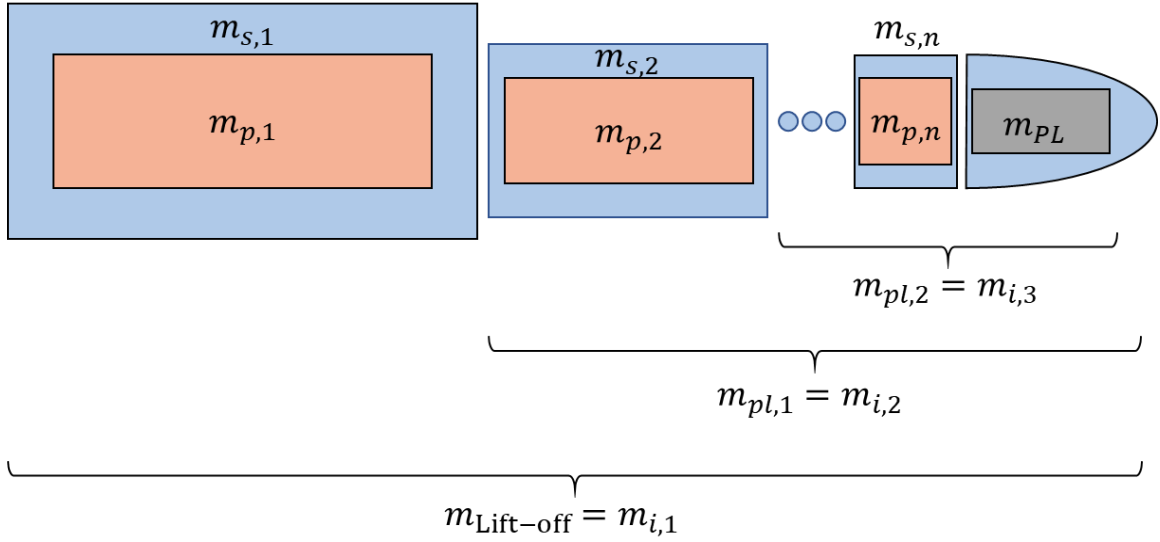

**Figure 1 Mass components of a multi-stage launch vehicle**

$$m_{pl,k} = m_{PL} + \Sigma_{j=k+1}^{n}(m_{s,j} + m_{p,j}) = \begin{cases} m_{i,k+1} & (k < n) \\ m_{PL} & (k = n) \end{cases} \qquad (2)$$

The mass ratio of each stage ($\mu_k$) is the ratio between the initial mass ($m_{i,k}$) and the final mass ($m_{f,k}$) after consuming the whole propellant of the stage.

$$\mu_k \equiv \frac{m_{i,k}}{m_{f,k}} = \frac{m_{i,k}}{m_{i,k} - m_{p,k}} = \frac{m_{s,k} + m_{p,k} + m_{pl,k}}{m_{s,k} + m_{pl,k}} \qquad (3)$$

Based on the definitions above, the relationship between the initial and payload masses of the $k^{th}$ stage can be derived as follows.

$$\frac{m_{i,k}}{m_{pl,k}} = \frac{m_{i,k}}{m_{i,k} - m_{p,k} - m_{s,k}} = \frac{\mu_k(1 - \varepsilon_k)}{1 - \mu_k \varepsilon_k} \qquad (4)$$

Therefore, the ratio of the lift-off mass to the payload mass ($m_{\text{Lift-off}}/m_{PL}$) is expressed as follows.

$$\frac{m_{\text{Lift-off}}}{m_{PL}} = \frac{m_{i,1}}{m_{pl,1}} \times \frac{m_{pl,1}}{m_{pl,2}} \times \cdots \times \frac{m_{pl,n}}{m_{PL}} = \frac{m_{i,1}}{m_{pl,1}} \times \frac{m_{i,2}}{m_{pl,2}} \times \cdots \times \frac{m_{i,n}}{m_{pl,n}} = \prod_{k=1}^{n} \frac{\mu_k(1 - \varepsilon_k)}{1 - \mu_k \varepsilon_k} \qquad (5)$$

The velocity increment generated by each stage can be calculated using Tsiolkovsky's rocket equation. The total velocity increment of the launch vehicle is obtained by summing the velocity increments produced by all stages as follows.

$$\Sigma_{k=1}^{n} \Delta v_k = \Sigma_{k=1}^{n} v_{ex,k} \ln \mu_k \qquad (6)$$

The velocity increment calculated by Eq. (6) should be greater than or equal to the required velocity increment ($\Delta v_{req}$) to put the payload into the target orbit. Hence, the optimal stage sizing problem can be formulated as follows.

[**P**$_S$: Optimal Stage Sizing]

For given values of $v_{ex,i}$, $\varepsilon_k$, and $\Delta v_{req}$

$$\max_{\mu_1,\cdots,\mu_n} \Pi_{k=1}^{n} \frac{\mu_k(1-\varepsilon_k)}{1-\mu_k\varepsilon_k} \tag{7}$$

subject to

$$\Sigma_{k=1}^{n} v_{ex,k} \ln \mu_k - \Delta v_{req} \geq 0. \tag{8}$$

One can solve the formulated optimization problem by introducing the Lagrange multiplier and applying a numerical root-finding procedure [10]. Based on the obtained optimal mass ratios ($\mu_k^*$) and the optimal mass and velocity increment of each stage ($m_{p,k}^*$, $m_{s,k}^*$, $\Delta v_k^*$) can be calculated.

## B. Considering Velocity Losses

While various studies have addressed the optimal stage sizing problem described in the previous section, early studies did not consider velocity losses. Instead, they determined the required velocity increment ($\Delta v_{req}$) using the target orbital velocity ($v_f$) and the initial velocity due to the Earth's rotation ($v_i$). However, the actual flight of a rocket involves the velocity losses caused by various sources (e.g., gravity, thrust vector control, aerodynamic drag, and pressure difference at the engine exit). The velocity loss reduces the final speed of the launch vehicle. Consequently, the velocity loss necessitates a higher velocity increment provided by the rocket, and its consideration for optimal stage sizing is critical.

Some studies assessed the required velocity increment ($\Delta v_{req}$) with $v_f$, $v_i$, and $\Delta v_{Loss}$ values ($\Delta v_{req} = v_f - v_i + \Delta v_{Loss}$) [10-11]. However, the velocity loss for each stage does not necessarily scale proportionately to the velocity increment allocated to the stage, potentially leading to a suboptimal design [12]. Jamilnia [10] and Jo [13] proposed a simultaneous optimization of staging and trajectory to consider the velocity losses along with staging. However, their work did not consider the realistic constraints for launch operations (e.g., impact points of separated stages and fairing) and lacked sufficient explanation/discussion on the implementation details of the simultaneous optimization approach.

More recent studies adopted a different approach [14-17]. Instead of simply adding the loss terms to $\Delta v_{req}$, they obtained the ideal velocity increment for each stage under the loss-free condition. Then they assigned the velocity loss from stage $k$ ($\Delta v_{Loss,k}$) to $\Delta v_k$. Since $\Delta v_{Loss,k}$ is directly linked to the flight trajectory, any modification to a stage size may lead to a change in velocity loss associated with the stage. Consequently, the methodology iteratively solves the trajectory optimization problem and the stage sizing problem, yielding an optimal launch vehicle size design. Although this methodology showed reasonable performance, the iterative process does not always obtain an optimal design due to the interdependencies among stage sizing, trajectory, and velocity losses. Particularly, one can observe that the process converges to a sub-optimal design when imposing constraints related to the launch safety operations.

Note that the new framework proposed in this paper can consider operational constraints simultaneously and efficiently obtain the optimal solution. The case study section of this paper provides the discussions on implementation details and the advantages of the proposed framework.

### C. Instantaneous Impact Point (IIP) for Launch Safety

An instantaneous impact point (IIP) is the touch-down point of a launch vehicle in case of the immediate termination of the propelled flight. Monitoring the IIP during launch operation is one of the most crucial tasks for ensuring the safety of space launch vehicles [18-19]. A multi-stage expendable launch vehicle separates the lower stage(s), which drops to the Earth's surface. The trajectory of a launch vehicle should be carefully designed so that the collision risk of the separated stage is minimized [20]. The position and velocity of a launch vehicle determine its impact point of a launch vehicle. Considering that the stage sizing determines the velocity increment, the integrated design of the vehicle stage and launch trajectory can effectively address the issue of operational safety along with the launch capability (i.e., maximized payload mass ratio) [21]. Figure 2 illustrates the concept of the IIP.

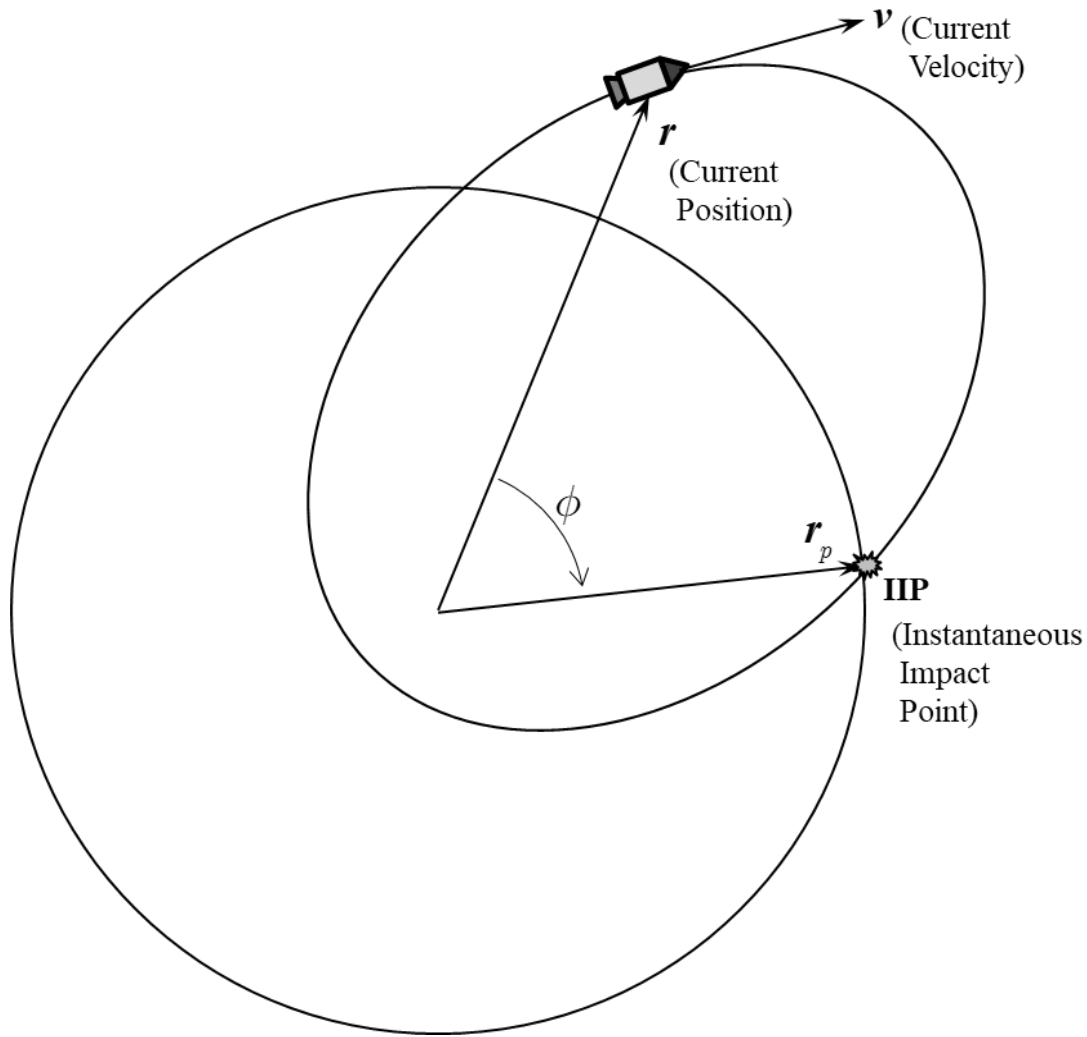

**Figure 2 Impact point of separated stage**

Accuracy and speed are two critical features required for IIP prediction used in a range safety system. Therefore, studies for IIP modeling and computation have focused on two primary objectives: enhancing the prediction accuracy and accelerating the computation. Montenbruck et al. [22-23] calculated the IIP based on the assumption that the launch vehicle experiences a constant downward gravitational acceleration in the LVLH coordinate system. They compensated for the modeling errors (e.g., earth curvature, gravitational acceleration, and aerodynamic drag) by introducing the correction terms for the calculation. On the other hand, the Keplerian algorithm used Kepler's two-body problem as the basis to predict the IIP. This methodology is considered more suitable for operating space launch vehicles than the past approach based on the constant-gravity assumption. The algorithm introduced in references [24] and [25] calculates the IIP while assuming that the vehicle's motion without propulsion is subject to only Kepler's central gravitational acceleration. These algorithms employed an iterative procedure to determine the IIP. Alternatively, Ahn and Roh [26] proposed a Keplerian algorithm without iterations, reducing the potential issues related to the computational burden caused by excessive iterations. Ahn and Seo [27] improved the algorithm to reflect the atmospheric drag and oblateness of the Earth based on the response surface method. The response surface method has been used in various aerospace applications to model physical quantities and explore the design space efficiently [28-30]. Their approach modeled the corrections terms for the flight angles and time with the current position and vector as inputs and improved the prediction accuracy of the IIP significantly. Nam and Ahn [31] introduced the concept of the adjusted instantaneous impact point (AIIP) to improve the flight safety (termination) decision-making process based on the IIP. The AIIP considers the potential delay between the flight termination command generation and the actual engine cut-off by using the "*Delta IIP*." They developed the analytical formula representing the time derivatives of IIP due to the external force (thrust) to compute the *Delta IIP*. They showed that AIIP-based decision-making could reduce the chances of undetected collision risk and unnecessary flight termination, improving flight safety operations.

The next section presents a framework that can optimize the launch vehicle stages and flight trajectory integratively while considering the range safety constraints on the IIP.

# III. Simultaneous Optimization Framework

The optimization framework proposed in this study consists of three modules. The first module is a flight simulation module that simulates the motion of a launch vehicle. It uses the specification of the launch vehicle, flight sequence, and the histories of attitude angles as input parameters and generates the vehicle's states (e.g., mass, position/velocity, and velocity loss terms) by integrating the dynamic equations. The second module takes the state variables (output of the first module) and computes output variables associated with flight operations (e.g., flight environment, range safety, orbital elements, and radar tracking). Various coordinate transformation and conversion formulae are used to calculate the output variables from the states. One of the important output variables associated with range safety operations is the instantaneous impact point (IIP), which is explained in detail in the following subsection. The last module formulates the simultaneous stage/trajectory design as a nonlinear optimization problem and solves it using a gradient-based algorithm.

The decision variables of the optimization problem are the structural/propellant masses of each stage ($m_{s,k}$ and $m_{p,k}$) and the angles representing the thrust directions at discretized time points ($\theta(t_j)$ and $\psi(t_j)$), which constitute the input variables of the first module. Its objective is to minimize the lift-off mass of the launch vehicle while satisfying the specified constraints on the target orbit insertion and launch operations. The constraints of the optimal design problem are imposed on the output variables of the second module (e.g., altitude, inertial speed, flightpath angle, inclination angle, and IIP) at discretized time points. The framework assumes the exhaust velocity ($v_{ex,k}$) and structural mass fraction ($\varepsilon_k$) as given constant parameters. Figure 2 shows the architecture of the simultaneous optimization framework proposed in this paper.

The proposed framework discretizes the entire flight time (starting from launch and ending with orbit insertion) into multiple intervals referred to as *Phases*. For example, Phase $p$ is the time interval starting/terminating at $t_{p-1}$ and $t_p$, respectively (= [$t_{p-1}$, $t_p$]). While trajectory optimization is an optimal control problem involving continuous states and input in nature, parameter discretization enables designers to handle the problem together with stage sizing in a single optimization framework. Note that the duration of a phase involving active thrust is subject to change depending on the propellant mass of each stage, while the duration of a coasting phase is constant. Important flight events (e.g., stage/fairing separations, ignition/cut-off of an engine) occur at the end of each phase ($t_p$). Figure 3 illustrates the configuration of a phase using discretized time points. The following subsections explain the three modules constituting the proposed framework.

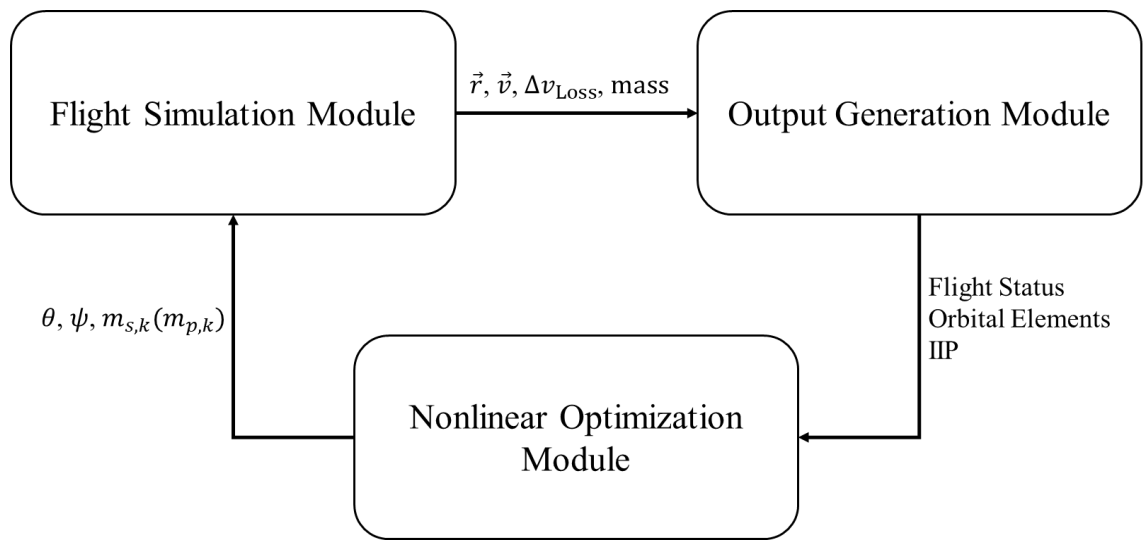

**Figure 3 Architecture of the proposed optimization framework**

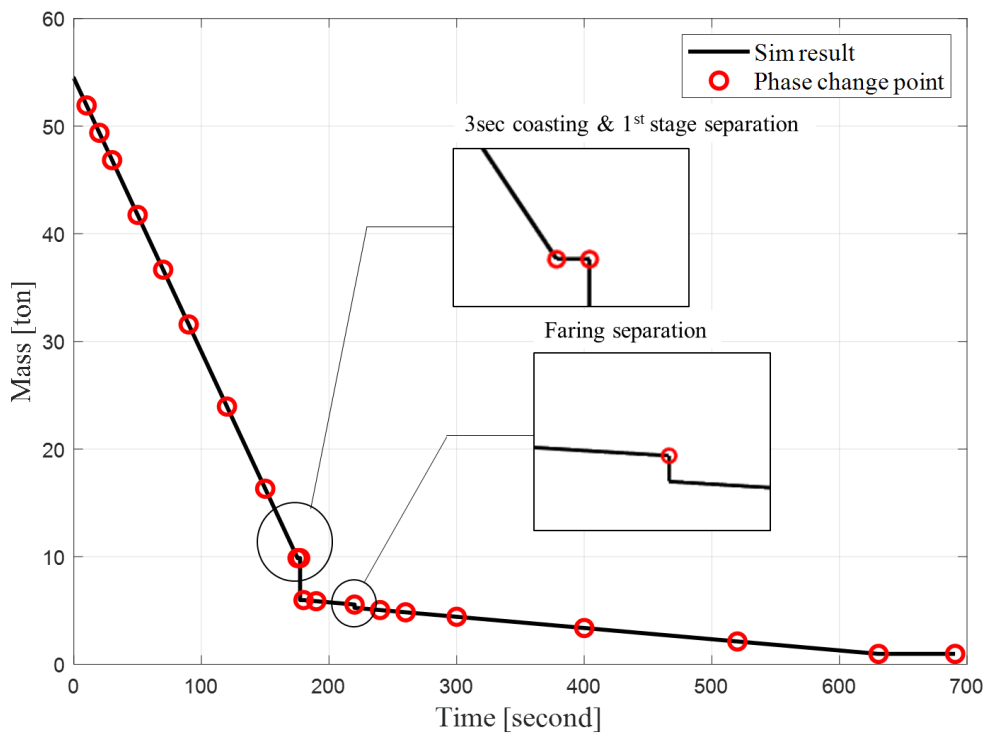

**Figure 4 Example of phase configuration**

## A. Flight Simulation Module

The flight simulation module implements a 3-degree-of-freedom point mass model. This module computes the position ($r$), velocity ($v$), and velocity loss values ($v_{Loss}$) over simulation time by integrating the differential equation describing the vehicle dynamics. The dynamic equations describing the vehicle motion at each phase are as follows.

$$\frac{d\boldsymbol{r}}{dt} = \boldsymbol{v} \tag{9}$$

$$\frac{d\boldsymbol{v}}{dt} = \left( \boldsymbol{f}_T + \boldsymbol{f}_a + \boldsymbol{f}_g \right) / m \tag{10}$$

$$\dot{m}_{p,k} = \frac{T_{v,k}}{g_0 I_{sp,k}} \tag{11}$$

$$\frac{d\Delta v_{L,p}}{dt} = \frac{p_a A_{e,j}}{m} \tag{12}$$

$$\frac{d\Delta v_{L,a}}{dt} = \frac{\| \boldsymbol{f}_a \|}{m} \tag{13}$$

$$\frac{d\Delta v_{L,g}}{dt} = \frac{\| \boldsymbol{f}_g \| \sin \gamma}{m} \tag{14}$$

$$\frac{d\Delta v_{L,TVC}}{dt} = \frac{\| \boldsymbol{f}_T \| (1 - \cos \delta)}{m} \tag{15}$$

Eqs. (9)-(11) represent the translational motion of the launch vehicle expressed in an Earth-Centered Inertial (ECI) coordinate frame. In these equations, m is the mass of the launch vehicle expressed as follows.

$$m = m_{PL} + \Sigma_{k=j}^{n} (m_{s,k} + m_{p,k}) + s \times m_{\text{faring}} \tag{16}$$

In Eq. (16), $j$ is the index representing the currently active stage, and $s$ is a binary indicator variable whose value is one before fairing separation and zero afterward. Three force components ($\boldsymbol{f}_a$, $\boldsymbol{f}_g$, and $\boldsymbol{f}_T$) are defined as follows.

$$\boldsymbol{f}_T = \left( v_{ex,j} \dot{m}_{p,j}(t) - A_{e,j} p_a \right) \boldsymbol{e}_T \tag{17}$$

$$f_a = -\frac{1}{2}\rho_a C_D S_{ref} \parallel \boldsymbol{v}_r \parallel \boldsymbol{v}_r \tag{18}$$

$$\boldsymbol{f}_g = m(\boldsymbol{g}_c + \boldsymbol{g}_p) = m\left\{-\frac{\mu}{\parallel \boldsymbol{r} \parallel^3}\boldsymbol{r} + \boldsymbol{g}_p\right\} \tag{19}$$

In these equations, $\boldsymbol{e}_T$ is a unit vector representing the thrust direction defined by the pitch and yaw angles ($\theta$ and $\psi$), and $p_a$ / $\rho_a$ are the atmospheric pressure/density from the 1976 U.S. standard atmosphere model [32]. Also, $\boldsymbol{g}_c$ and $\boldsymbol{g}_p$ are the central and perturbed (considering $J_2$ term) gravity acceleration vectors, respectively. Eqs. (12)-(15) are the differential equations for velocity loss components: pressure loss ($\Delta v_{L,p}$), drag loss ($\Delta v_{L,a}$), gravity loss ($\Delta v_{L,g}$), and TVC loss ($\Delta v_{L,TVC}$), respectively.

The launch position and velocity determine the initial states of the first phase. Also, the final states of phase p determine the initial states of phase ($p$+1). These relationships are expressed as follows.

$$\begin{bmatrix} \boldsymbol{r}_i^{(p+1)} \\ \boldsymbol{v}_i^{(p+1)} \end{bmatrix} = \begin{bmatrix} \boldsymbol{r}_f^{(p)} \\ \boldsymbol{v}_f^{(p)} \end{bmatrix} = \begin{bmatrix} \boldsymbol{r}(t_p) \\ \boldsymbol{v}(t_p) \end{bmatrix} \tag{20}$$

## B. Output Generation Module

The output generation module computes various output variables using the state variables ($\boldsymbol{r}$, $\boldsymbol{v}$, $dm_k/dt$), which are the outputs of the first module. The module first computes the position and velocity elements expressed in various coordinate frames using coordinate transformations. The orbital elements are calculated by transforming the inertial position/velocity into six parameters [33]. The flight-safety-related variables (e.g., IIP and IIP rate) are computed based on the procedures presented in references [26]. Table 1 outlines the selected output variables (notations, descriptions, and units) computed in the second module.

**Table 1 Selected output variables from the second module**

| Category | Output | Description / Unit |
|---|---|---|
| **Flight Status** | $v_I$ | Inertial speed of the vehicle, *m/s* |
| | $v_R$ | Relative (to the Earth) speed of the vehicle, *m/s* |
| | $h$ | Altitude, *m* |
| | $\lambda$ | Longitude, *rad* |
| | $\varphi$ | Latitude, *rad* |
| | $q$ | Dynamic Pressure, *Pa* |
| | $\gamma$ | Flight path angle, *rad* |
| | $\chi$ | Azimuth angle, *rad* |
| | $\alpha$ | Angle of attack, *rad* |
| | $\beta$ | Side slip angle, *rad* |
| **IIP** | $t_{go}$ | Remaining time to flight until the IIP, *s* |
| | $\varphi_I$ | Latitude of IIP, *rad* |
| | $\lambda_I$ | Longitude of IIP, *rad* |
| **Orbital Element** | $i$ | Inclination angle, *rad* |
| | $h_p$ | Perigee altitude, *m* |
| | $h_a$ | Apogee altitude, *m* |
| | $a$ | Semi-major axis, *m* |
| | $f$ | True anomaly, *rad* |

## C. Nonlinear Optimization Module

The nonlinear optimization module formulates and optimizes the problem of minimizing lift-off mass with constraints on the target orbit insertion and other launch operations (e.g., flight safety). The problem is formulated as follows.

[**P**$_I$: Integrated Optimization of Stage and Vehicle]

$$\min_{x} m_{\text{Lift-off}} \tag{21}$$

$$\boldsymbol{x} = [m_{s,k}, \theta_i^p, \psi_i^p] \quad (1 \le k \le N, \ 1 \le p \le P) \tag{22}$$

subject to

$$\varepsilon_k = m_{s,k} \, / \, (m_{s,k} + m_{p,k}) \tag{23}$$

Eqs. (9)-(20)

$$\begin{bmatrix} h_f^P \\ \gamma_f^P \\ (v_I)_f^P \\ i_f^P \end{bmatrix} = \begin{bmatrix} h_{req} \\ \gamma_{req} \\ (v_I)_{req} \\ i_{req} \end{bmatrix} \tag{24}$$

$$\alpha^p = 0, \ (p \in \Pi_G) \tag{25}$$

$$q(t) \le q_{max} \tag{26}$$

$$(\varphi_I^k)_L \le \varphi_I^{S_k} \le (\varphi_I^k)_U, (\varphi_I^F)_L \le \varphi_I^{S_F} \le (\varphi_I^F)_U \tag{27}$$

$$(\lambda_I^k)_L \le \lambda_I^{S_k} \le (\lambda_I^k)_U, (\lambda_I^F)_L \le \lambda_I^{S_F} \le (\lambda_I^F)_U \tag{28}$$

Eq. (21) is the objective function of the nonlinear optimization problem, which is the minimization of the lift-off mass of the vehicle for a given target mission and payload mass. Eq. (22) defines the decision variables for the problem: structural mass of stages ($m_{s,k}$) and discretized angles defining thrust direction ($\theta_f^p$ and $\psi_f^p$). We select the final pitch and yaw angles at the end of each phase for design variables and use linear interpolation to determine the angles within the phase. Eq. (23) represents a constraint on the structural mass fraction ($\varepsilon_k$), determining the propellant mass ($m_{p,k}$). Vehicle dynamics expressed as Eqs. (9)-(20) are imposed as constraints as well. The differential equations are integrated numerically to create the output variables representing the orbital elements and launch operations, which are used to define constraints expressed in Eqs. (24)-(28). Eq. (24) includes the constraints for target orbit insertion specifying the final altitude, flightpath angle, inertial velocity, and inclination angle. Eq. (25) represents the constraint that the angle of attack should be zero during the gravity turn, where $\Pi_G$ is the set of phases that belong to the gravity turn period. Eq. (26) describes the constraint on the dynamic pressure the vehicle experiences throughout the ascending flight. Eqs. (27)-(28) impose the constraints on the impact point(s) at the end of separation stages by their lower and upper bounds. In these equations, superscript $S_k$ and $S_F$ represent the phases at which the $k^{th}$ stage and payload fairing are separated, respectively.

Sometimes the designer may impose additional constraints on launch operations, such as the structural load (expressed by the product of the dynamics pressure and aerodynamic angle) and thermal heating. Equalities or inequalities (e.g., $\boldsymbol{g(x)} \le \boldsymbol{0}, \boldsymbol{h(x)} = \boldsymbol{0}$) can express these additional constraints.

## IV. Case Study

This section presents a case study demonstrating the effectiveness of the proposed framework to optimize the stage sizing and trajectory of a launch vehicle simultaneously. We set up the mission of delivering a 3000 kg payload to a circular target orbit (altitude: 300 km, inclination: 80 deg). The characteristics of the structural and propulsion systems for Korean Space Launch Vehicle II (KSLV-II, Nuri) are used for this case study (e.g., structural mass fraction, specific impulse, and thrust). Table 2 summarizes the specifications of the KSLV-II vehicle used for this case study. The structural mass values presented in this table were used as the initial guess for the decision variables.

**Table 2 Reference design for case study - specification of KSLV-II**

| Parameter | Stage 1 | Stage 2 | Stage 3 | Payload | Faring |
|---|---|---|---|---|---|
| Thrust, $tonf$ | 303.2 | 75.8 | 7 | | |
| $v_{ex}$, $m/s$ | 2923 | 3093 | 3188 | | |
| $\dot{m}$, $kg/s$ | 1017 | 240 | 21.5 | | |
| $A_e$, $m^2$ | 3.6 | 0.9 | 0.2 | | |
| $m_s$, $kg$ | 14900 | 5300 | 1800 | | |
| $m_p$, $kg$ | 128200 | 36600 | 10800 | | |
| Stage mass, $kg$ | 143100 | 41900 | 12600 | 3000 | 900 |
| $\varepsilon$, - | 0.104 | 0.127 | 0.143 | | |
| $S_{ref}$, $m^2$ | 9.6 | 5.3 | 5.3 | | |
| $C_D$, - | $0.35 - 1.1$ | $0.35 - 1.1$ | - | | |
| Launch Site | | 127.5 °E / 34.4 °N / 140 $m$ | | | |

We implemented the proposed simultaneous optimization framework and the method proposed by Koch [14] for comparison. Koch's algorithm requires iterations between trajectory optimization and staging problems. The convergence criterion for Koch's algorithm was selected as 0.1 % error for the maximum payload obtained from the trajectory optimization. The case study uses MATLAB on a Windows 11 operating system with an Intel i5-13600K processor and 64 GB of RAM. The sequential quadratic programming (SQP) provided by MATLAB Optimization Toolbox was used as the algorithm used for the optimization. We conducted two case studies on redesigning KSLV-II using the proposed framework. Case I considers the gravity turn condition for the operational constraint (Eq. (25)), Case II additionally considers the dynamic pressure constraint (Eqs. (25)-(26)), and Case III considers the gravity turn and the IIP range constraints (Eqs. (25), (27), and (28)) for the optimization.

Table 3 exhibits the optimization results for Case I obtained using the traditional stage-trajectory iteration and the simultaneous optimization. The simultaneous optimization method presented in this paper showed a distinct advantage by treating both the stage mass and attitude as a design variable within a single nonlinear optimization problem. This approach leads to improved solutions with a lighter objective function compared to a traditional approach that handles staging and trajectory optimization separately. The simultaneous optimization method resulted in a lower $\Delta v_{Loss}$ and a 7.2% reduction in the lift-off mass than the traditional algorithm, supporting the advantage of the proposed framework.

**Table 3 Case I (Constraint on gravity turn) results**

| Parameter | Staging – Trajectory Iteration | Simultaneous Optimization |
|---|---|---|
| $V_f$ -$V_i$, $m/s$ | 7342 | 7342 |
| $\Delta v_{Loss}$, $m/s$ | 2019 | 1951 |
| $\Delta v_k$, $m/s$ | [3871, 2712, 2777] | [3296, 2860, 3136] |
| $m_{s,k}$, $kg$ | [14,174, 2498, 908] | [12,116, 3263, 1161] |
| $m_{p,i}$, $kg$ | [121,950, 17,249, 5448] | [104,246, 22,531, 6968] |
| Lift-off mass, $kg$ | 166,161 | 154,185 |
| Number of iterations, - | 11 | - |
| Run time, $s$ | 90 | 31 |

Table 4 presents the results for Case II. The maximum dynamic pressure value used for the constraint ($q_{max}$) is 40 $kPa$. One can observe that the proposed framework yields the objective function (lift-off mass) 6.8 % smaller than the traditional approach involving the iteration between the trajectory design and staging. In addition, the computation time spent on the proposed framework (28 $s$) was significantly shorter than that of the traditional approach (4766 $s$), which is caused by the large number of iterations (405). Figure 5 shows the optimal trajectories obtained by both methodologies. The launch vehicle design obtained through the simultaneous optimization framework moves slightly more vertically.

**Table 4 Case II (Constraints on gravity turn and dynamic pressure) results**

| Parameter | Staging – Trajectory Iteration | Simultaneous Optimization |
|---|---|---|
| $V_f$ -$V_i$, $m/s$ | 7342 | 7342 |
| $\Delta v_{Loss}$, $m/s$ | 2156 | 2114 |
| $\Delta v_k$, $m/s$ | [3980, 2815, 2703] | [3073, 3264, 3119] |
| $m_{s,k}$, $kg$ | [15,417, 2580, 858] | [12,567, 4251, 1147] |
| $m_{p,i}$, $kg$ | [132,648, 17,815, 5146] | [108,125, 29,353, 6882] |
| Lift-off mass, $kg$ | 178,362 | 166,224 |
| Number of iterations, - | 405 | - |
| Run time, $s$ | 4766 | 28 |

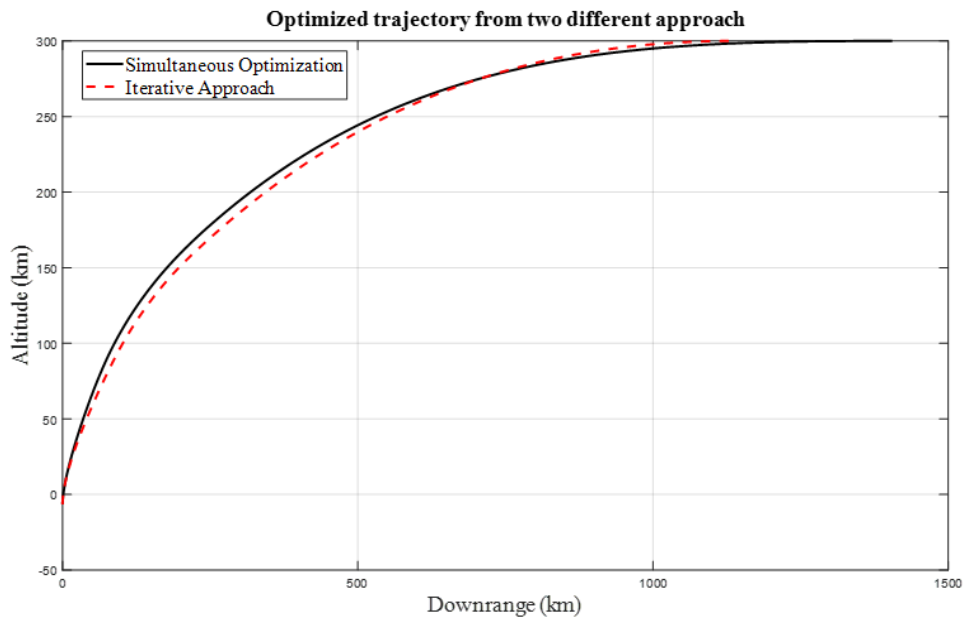

**Figure 5 Optimized trajectories obtained by the two methods**

In the conventional iterative approach, the value of velocity losses can vary significantly as the stage size of the launch vehicle changes during iteration. Reduction in the lift-off mass results in significant acceleration during the early phases after launch, leading to a relatively higher speed at the same altitude and an increase in the dynamics pressure. In order to satisfy the dynamic pressure constraint, the flight trajectory needs to be adjusted, varying the velocity loss terms. Thus, the operational constraints serve as incentives for changes at each stage mass during the design iteration process. The magnitude of these changes can significantly impact the convergence time of the traditional approach or even fail to converge in the worst-case scenario. In contrast, the proposed framework solves the nonlinear optimization problem by simultaneously considering the trajectory and staging. Therefore, the framework can successfully reflect operational constraints in staging decisions, yielding the solution more efficiently.

Table 5 presents the results for Case III, which constrains the IIP of the separated first stage. This constraint ensures that the separated tank falls within the safe area and does not pose risks to human populations or cross borders. The upper limit of the IIP longitude is set as 128.3 °E, and the zero angle of attack constraint during the gravity turn period is imposed as well.

**Table 5 Case III (Constraints on gravity turn and IIP of Stage 1) results**

| Parameter | Staging – Trajectory Iteration | Simultaneous Optimization |
|---|---|---|
| $V_f$ -$V_i$, $m/s$ | - | 7342 |
| $\Delta v_{Loss}$, $m/s$ | - | 1977 |
| $\Delta v_k$, $m/s$ | - | [3076, 3089, 3154] |
| $m_{s,k}$, $kg$ | - | [11,755, 3841, 1175] |
| $m_{p,i}$, $kg$ | - | [101,141, 26,525, 7050] |
| Lift-off mass, $kg$ | - | 155,388 |
| Number of iterations, - | Diverged after Iteration 6 | - |
| Run time, $s$ | - | 34.32 |

Figure 6 visualizes the trajectory of the current position and IIP of the vehicle on the map. Since the target orbit is prograde, the launch vehicle flies eastward. Therefore, as the first stage size increases, the separated first stage's impact point shifts eastward. Hence, when there is a constraint on the maximum allowable longitude of the first stage impact point, it is necessary to design the first stage smaller so that the velocity increment allocated to the first stage decreases. The velocity increment values achieved by the first stage are 3296 m/s (Case I) and 3076 m/s (Case III), respectively, showing this tendency. Note that the traditional optimization approach failed to obtain the solution for Case III (procedure diverged after iteration 6), showing its vulnerability to design problems involving complex operational constraints.

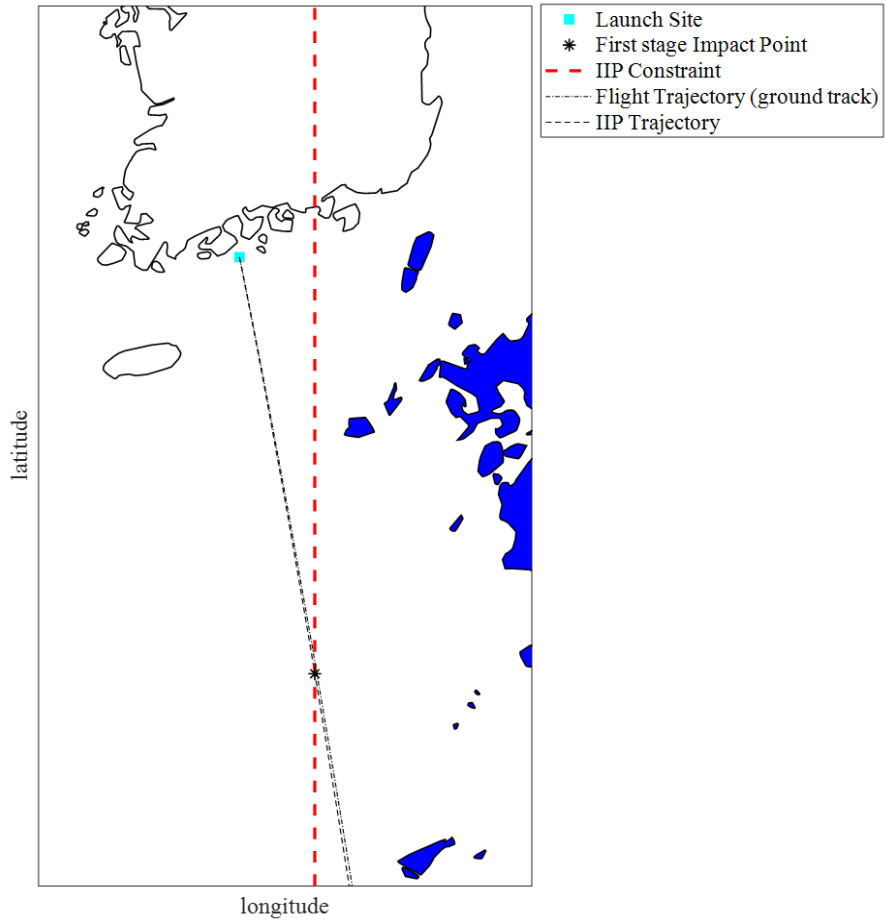

**Figure 6 Trajectories (position and IIP) and IIP constraint shown on the map**

## V. Conclusions

This study proposed a framework to optimize a launch vehicle's trajectory and stage sizing simultaneously. In particular, the proposed simultaneous optimization framework can effectively handle the constraints on the launch operations of the vehicle. The framework can formulate the integrated design problem as a single nonlinear optimization problem by properly discretizing the control variables for ascending flight and combining them with the optimal sizing problem. Unlike the iterative approach between trajectory optimization and stage sizing adopted by previous studies, the proposed framework solves a single nonlinear optimization problem, improving its convergence characteristics and computational speed. Three case studies support the advantages of the proposed framework mentioned above. An extension of the proposed framework to address design problems involving reusable launch vehicles is a potential subject for future study, which requires the modification of the trajectory simulation and optimization modules to address the returning/landing trajectory of the reusable stage.

# References

[1]  Balesdent, M., Bérend, N.,  Dépincé, P., and Chriette, A., "A survey of multidisciplinary design optimization methods in launch vehicle design," *Structural and Multidisciplinary Optimization*, Vol. 45, No. 5., 2012, pp. 619–642.

https://doi.org/10.1007/s00158-011-0701-4

[2]  Shu, J. I., Kim, J. W., Lee, J. W., and Kim, S., "Multidisciplinary mission design optimization for space launch vehicles based on sequential design process," *Proceedings of the Institution of Mechanical Engineers, Parg G: Journal of Aerospace Engineering*, Vol. 230, No. 1, 2016, pp. 3–18.

https://doi.org/10.1177/0954410015586858

[3]  Malina, F., and Summerfield, M., "The Problem of Escape from the Earth by Rocket," *Journal of the Aeronautical Sciences*, Vol. 14, No. 8, 1947, pp. 471–480.

https://doi.org/10.2514/8.1417

[4]  Goldsmith, M., "On the Optimization of Two-Stage Rockets," *Journal of Jet Propulsion*, Vol. 27, No. 4, 1957, pp. 415–428.

https://doi.org/10.2514/8.12799

[5]  Schurmann, E., "Optimum Staging Technique for Multistaged Rocket Vehicles," *Journal of Jet Propulsion*, Vol. 27, No. 8, 1957, pp. 863–865.

https://doi.org/10.2514/8.12965

[6]  Hall. H., and Zambelli, E., "On the Optimization of Multi-stage Rockets," *Journal of Jet Propulsion*, Vol. 28, No. 7, 1958, pp. 463–465.

https://doi.org/10.2514/8.7353

[7]  Subotowicz, M., "The Optimization of the N-Step Rocket With Different Construction Parameters and Propellant Specific Impulses in Each Stage," Journal of Jet Propulsion, Vol. 28, No. 7, 1958, pp. 460–463.

https://doi.org/10.2514/8.7352

[8]  Ragsac, R., and Patterson, P., "Multi-stage Rocket Optimization," *ARS Journal*, Vol. 31, No. 3, 1961, pp. 450–452.

https://doi.org/10.2514/8.5509

[9]  Fan, L. T., and Wan, C. G., "Weight minimization of a step rocket by the discrete maximum principle," *Journal of Spacecraft and Rockets*, Vol. 1, No. 1, 1964, pp. 123–125.

https://doi.org/10.2514/3.27608

[10]  Jamilnia, R., and Naghash, A., "Simultaneous optimization of staging and trajectory of launch vehicles using two different approaches," *Aerospace Science and Technology*, Vol. 23, No. 1, 2012, pp. 85–92.

https://doi.org/10.1016/j.ast.2011.06.013

[11]  Civek-Coşkun, E., and Özgören, K., "A generalized staging optimization program for space launch vehicles," in *2013 6th International Conference on Recent Advances in Space Technologies* (*RAST*), 2013, pp. 857–862.

https://doi.org/10.1109/RAST.2013.6581333

[12]  Edberg, D., and Costa, W., *Design of Rockets and Space Launch Vehicles*, Second Edition, American Institute of Aeronautics and Astronautics, Reston, VA, USA, 2022, pp. 271–336.

https://doi.org/10.2514/5.9781624106422.0271.0336

[13]  Jo, B.U., and Ho, K., "Simultaneous Sizing of a Rocket Family with Embedded Trajectory Optimization," ArXiv Preprint ArXiv:2302.12900, 2023.

https://doi.org/10.48550/arXiv.2302.12900

[14]  Koch, A. D., "Optimal staging of serially staged rockets with velocity losses and fairing separation," *Aerospace Science and Technology*, Vol. 88, 2019, pp. 65–72.

https://doi.org/10.1016/j.ast.2019.03.019

[15]  Jo, B. U., and Ahn, J., "Optimal staging of reusable launch vehicles considering velocity losses," *Aerospace Science and Technology*, Vol. 109, 2021, 106431.

https://doi.org/10.1016/j.ast.2020.106431

[16]  Jo, B. U., and Ahn, J., "Optimal staging of reusable launch vehicles for minimum life cycle cost," *Aerospace Science and Technology*, Vol. 127, 2022, 107703.

https://doi.org/10.1016/j.ast.2022.107703

[17]  Cho, B., Jo, B. U., and Ahn, J., "Integrated Framework for Staging and Trajectory Optimization of a Launch Vehicle Considering Range Safety Operations," *International Journal of Aeronautical and Space Sciences*, Vol. 22, No. 4, 2021, pp. 963–973.

https://doi.org/10.1007/s42405-020-00348-6

[18]  Eastern and Western Range 127-1 Range Safety Requirements, U.S. Air Force, Patrick AFB Range Safety Office, Brevard, FL, 1999.

[19]  Committee on Space Launch Range Safety, Aeronautics and Space Engineering Board, National Research Council, Streamlining Space Launch Range Safety, National Academic Press, Washington, D.C., 2000, p. 19.

[20]  Yoon, N., and Ahn, J., "Trajectory Optimization of a Launch Vehicle with Explicit Instantaneous Impact Point Constraints for Various Range Safety Requirements," Journal of Aerospace Engineering, Vol. 29, No. 3, 2016, 06015003.

https://doi.org/10.1061/(ASCE)AS.1943-5525.0000567

[21]  Cho, B., Jo, B. U., and Ahn, J., "Integrated Framework for Staging and Trajectory Optimization of a Launch Vehicle Considering Range Safety Operations," International Journal of Aeronautical and Space Sciences, Vol. 22, No. 4, 2021, pp. 963–973.

https://doi.org/10.1007/s42405-020-00348-6.

[22]  Montenbruck, O., and Markgraf, M., "Global positioning system sensor with instantaneous-impact-point prediction for sounding rockets," *Journal of Spacecraft and Rockets*, Vol. 41, No. 4, 2004, pp. 644–650.

https://doi.org/10.2514/1.1962

[23] Montenbruck, O., Markgraf, M., Jung, W., Bull, B., and Engler, W., "GPS based prediction of the instantaneous impact point for sounding rockets," *Aerospace Science and Technology*, Vol. 6, No. 4, 2002, pp. 283-294.

https://doi.org/10.1016/S1270-9638(02)01163-X

[24] Licensing and Safety Requirements for Operation of a Launch Site, Federal Aviation Admin. Code of Federal Regulations 14, Washington, D.C., Oct. 2000, Parts 401, 417, 420.

[25] Program to Optimize Simulated Trajectories (POST): Formulation Manual, Vol. 1, Martin Marietta Corp., Baltimore, MD, April 1975.

[26] Ahn, J., and Roh, W. -R., "Noniterative instantaneous impact point prediction algorithm for launch operations," *Journal of Guidance, Control, and Dynamics*, Vol. 35, No. 2, 2012, pp. 645–648.

https://doi.org/10.2514/1.56395

[27] Ahn, J., and Seo, J., "Instantaneous Impact Point Prediction Using the Response-Surface Method," *Journal of Guidance, Control, and Dynamics*, Vol. 36, No. 4, 2013, pp. 958–966.

https://doi.org/10.2514/1.59625

[28] Kim, J., and Ahn, J., "Modeling and optimization of a reluctance accelerator using DOE-based response surface methodology," *Journal of Mechanical Science and Technology*, Vol. 31, 2017, pp. 1321-1330.

https://doi.org/10.1007/s12206-017-0231-0

[29] Choi, U., Sung, T., and Ahn, J., "Adaptive Experimental Design for Aerodynamic Modeling with Hard-to-Change Factors," *Journal of Aerospace Information Systems*, Vol. 17, No. 6, 2020, pp. 311-315.

https://doi.org/10.2514/1.I010813

[30] Lee, M., and Ahn, J., "Flush Air Data System Modeling Using DOE-Based Wind-Tunnel Test," *International Journal of Aeronautical and Space Sciences*, Vol. 24, 2023, pp. 395-410.

https://doi.org/10.1007/s42405-022-00549-1

[31] Nam, Y., Seong, T., and Ahn, J., "Adjusted Instantaneous Impact Point and New Flight Safety Decision Rule," *Journal of Spacecraft and Rockets*, Vol. 53, No. 4, 2016, pp. 766-773.

https://doi.org/10.2514/1.A33424

[32]  *U. S. Standard Atmosphere, 1976*, NASA-TM-X-74335, National Aeronautics and Space Administration (NASA), US Government Printing Office, Washington, DC, 1976.

[33]  Vallado, D. A., *Fundamentals of Astrodynamics and Applications*, 4[th] ed., Microcosm Press, Portland, OR, 2013.

[34]  Ahn, J., and Roh, W., "Analytic Time Derivatives of Instantaneous Impact Point," *Journal of Guidance, Control, and Dynamics*, Vol. 37, No. 2, 2014, pp. 383-390.

https://doi.org/10.2514/1.61681

[35]  Jo, B., and Ahn, J., "New Formulation for Time Derivatives of Instantaneous Impact Point Based on Geometric Decomposition," *Journal of Aerospace Engineering*, Vol. 31, No. 4, July 2018, 04018036.

https://doi.org/10.1061/(ASCE)AS.1943-5525.0000865

[36]  Roh, W., Cho, S., Sun, B., Choi, K., Jung, D., Park, C., Oh, J., and Park, T., "Mission and System Design Status of Korea Space Launch Vehicle-II succeeding Naro Launch Vehicle," in Proceedings of the KSAS Fall Conference, 2012, pp. 233–239.